\documentclass{article}
\usepackage[english]{babel}

\usepackage{amsmath}
\usepackage{amssymb}
\usepackage{amsfonts}
\usepackage{mathrsfs}
\usepackage{graphics}

\newcommand{\ebox}{\fbox {} \smallskip}

\def\{{\protect\lbrace}
\def\}{\protect\rbrace}

\def\sign{\operatorname{sign}}

\newcommand{\pa}{\partial}

\def\alf{\alpha}
\def\bet{\beta}

\def\eps{\varepsilon}
\def\gam{\gamma}
\def\La{\Lambda}
\def\la{\lambda}
\def\phi{\varphi}

\def\ti{\tilde}

\def\tphi{\ti \phi}
\def\tpsi{\ti \psi}

\newcommand{\bR}{\mathbb{R}}

\newcommand{\bZ}{\mathbb{Z}}

\newcommand{\ee}{e}
\renewcommand{\O}{O}

\newcommand{\FF}{\mathscr {F}}
\newcommand{\KK}{\mathscr {K}}
\newcommand{\LL}{\mathscr {L}}
\newcommand{\CC}{{\gam_{1-b}}}
\newcommand{\CCC}{{\gam_{b}}}

\title{On pleated singular points of first order implicit differential equations}

\author{
R.A.~Chertovskih
\thanks{Centre for Wind Energy and Atmospheric Flows, Faculdade de Engenharia da Universidade do Porto, Rua Dr. Roberto Frias s/n, 4200-465 Porto, Portugal,
{\tt roman@mitp.ru}}
\thanks{International Institute of Earthquake Prediction Theory and Mathematical Geophysics, Profsoyuznaya str. 84/32, 117997  Moscow, Russia},
A.O.~Remizov
\thanks{CMAP \'Ecole Polytechnique Route de Saclay, 91128 Palaiseau Cedex, France,
{\tt alexey-remizov@yandex.ru}}
}

\begin{document}
\large
\maketitle

\begin{abstract}
We study phase portraits of a first order implicit differential equation in a neighborhood of its pleated singular point that is a non-degenerate singular point of the lifted field. Although there is no a visible local classification of implicit differential equations at pleated singular points (even in the topological category),
we show that there exist only six essentially different phase portraits, which are presented.
\end{abstract}

\section{Introduction}

A well-known geometrical approach to study implicit differential equations
\begin{equation}
F(x,y,p) = 0,   \ \ \ p = \frac{dy}{dx},
\label{1}
\end{equation}
consists of the lift the multivalued direction field defined by equation \eqref{1} on the $(x,y)$-plane to a single-valued direction field $X$ (which is called the lifted field) on the surface $\FF$ given by the equation $F(x,y,p) = 0$ in the $(x,y,p)$-space.\footnote{ This approach goes back to H.\,Poincar\'e (M\'emoire sur les courbes d\'efinies par les \'equations diff\'erentielles. -- J.~Math. Pures Appl., S\'er~4, 1885) and A.\,Clebsch (Ueber eine Fundamentalaufgabe der Invariantentheorie. --  G\"ottingen Abh.~XVII, 1872).
The latter paper contains geometric interpretation of differential equations (both ordinary and partial) and related theory of connexes, which is quite similar to the lifting; an account of these ideas is contained also in the famous book ``Vorlesungen \"uber h\"ohere Geometrie" by F.\,Klein.}
In this paper, the function $F$ is always supposed to be $C^{\infty}$-smooth. The lifted field $X$ is an intersection of the contact planes $dy = pdx$ with the tangent planes to the surface $\FF$, that is, $X$ is defined by the vector field
\begin{equation}
\dot x = F_p, \ \ \  \dot y = pF_p, \ \ \ \dot p = -(F_x+pF_y),
\label{2}
\end{equation}
whose integral curves are 1-graphs of integral curves (briefly, solutions) of equation \eqref{1}.
Conversely, solutions of \eqref{1} are $\pi$-projections of integral curves of $X$,
where $\pi$ is the projection from the surface $\FF$ to the $(x,y)$-plane along the $p$-direction (called vertical).

\begin{figure}[h!]
\begin{center}
\includegraphics{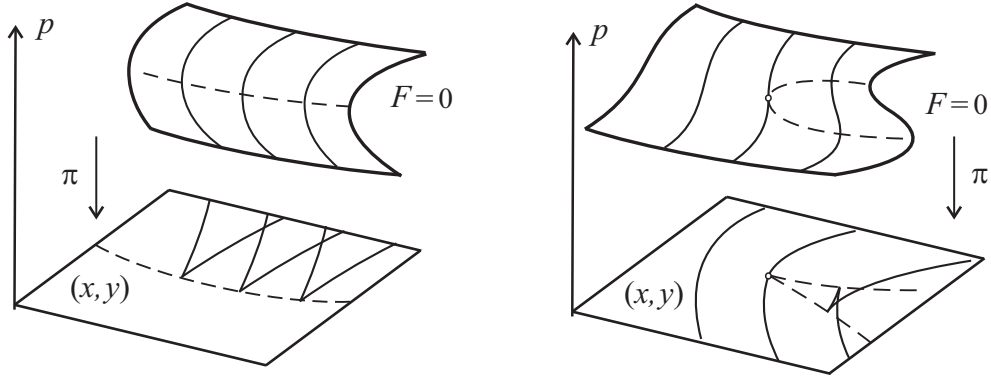}
\caption{
The projection $\pi$ has a fold (left) and pleat (right); the lifted field $X$ is defined on the surface $\FF$.
The dashed lines are the criminant and the discriminant curve on $\FF$ and the $(x,y)$-plane, respectively.
}\label{Fig00}
\end{center}
\end{figure}

This approach can be used for studying the local behavior of solutions of \eqref{1} near so-called \textit{singular points} -- points of the surface $\FF$ where $F_p=0$, that is, equation \eqref{1} cannot be locally resolved with respect to $p$ by the implicit function theorem and the germ of $\pi$ is not a diffeomorphism; see Fig.~\ref{Fig00}.

The first results of this sort were obtained by H.\,Poincar\'e, and later on, by various authors.
Moreover, this method allows to get a list of local normal forms of equation \eqref{1}.
Recall that two implicit differential equations are called smoothly (topologically) equivalent if there exists a diffeomorphism (homeomorphism, respectively) of the $(x,y)$-plane that sends integral curves of the first equation to integral curves of the second one.

\medskip

To describe the main results of this sort, we need to give several definitions following
\cite{A-Geom, A-Contact, AI, Dav-Japan, Rem}.

The locus of singular points of equation \eqref{1} is called the \textit{criminant} $\KK$; that is, $\KK$ is a set given by the equalities $F=F_p=0$.
The projection $\pi (\KK)$ on the $(x,y)$-plane is called the \textit{discriminant curve}.
The set $\LL$ given by the equalities $F=0 = F_x+pF_y=0$ is called the \textit{inflection curve}.\footnote{
The meaning of this name is clear from what follows.
Let $\gam$ be an integral curve of $X$, that is, an integral curve of the vector field \eqref{2}.
Suppose that the corresponding solution $\pi(\gam)$ of equation \eqref{1} has an inflection at some point on the $(x,y)$-plane.
Then the last component, $-(F_x+pF_y)$, of the vector field \eqref{2} vanishes at the corresponding point of the surface $\FF$.
}

We call a singular point $\O$ \textit{proper} if $F_x+pF_y \neq 0$, that is, the third component of the vector field \eqref{2} does not vanish at $\O$. Otherwise, we call it \textit{improper} singular point.

Improper singular points belong to the intersection $\KK \cap \LL$; they are characterized by the condition that
the surface $\FF$ is not regular or it is tangent to the contact plane.
Without loss of generality further we always assume $\O$ to be the origin in the $(x,y,p)$-space
(this can be obtained by appropriate affine map of the $(x,y)$-plane).

\medskip

A generic germ of equation \eqref{1} has singular points of the following three types.

1. \textit{A folded proper point}\footnote{In \cite{A-Geom, AI, A-Contact} such points are called \textit{regular} although being  singular points of implicit differential equation. However, we prefer to use another terminology.}:
the conditions $F_{pp} \neq 0$ and $F_x \neq 0$ hold true at $\O$.
Then the lifted field $X$ is defined, the criminant $\KK$ is regular and not vertical at $\O$, and the projection $\pi$ has a fold at all points of $\KK$. In a neighborhood of $\O$ each integral curve of $X$ transversally intersects $\KK$, and the corresponding solution of equation \eqref{1} has a cusp on the discriminant curve; see Fig.~\ref{Fig00} (left). Moreover, the whole family of solutions of \eqref{1} can be brought to the normal form $(y-c)^2=x^3$, $c \in \bR$, by a $C^{\infty}$-smooth diffeomorphism of the $(x,y)$-plane preserving the point $\O$.
The corresponding normal form $p^2 = x$ of equation \eqref{1} is named after Italian mathematician Maria Cibrario
(who established it in the analytic category) \cite{A-Geom, AI, Cib, Dav-85}.

\begin{figure}[h!]
\begin{center}
\includegraphics{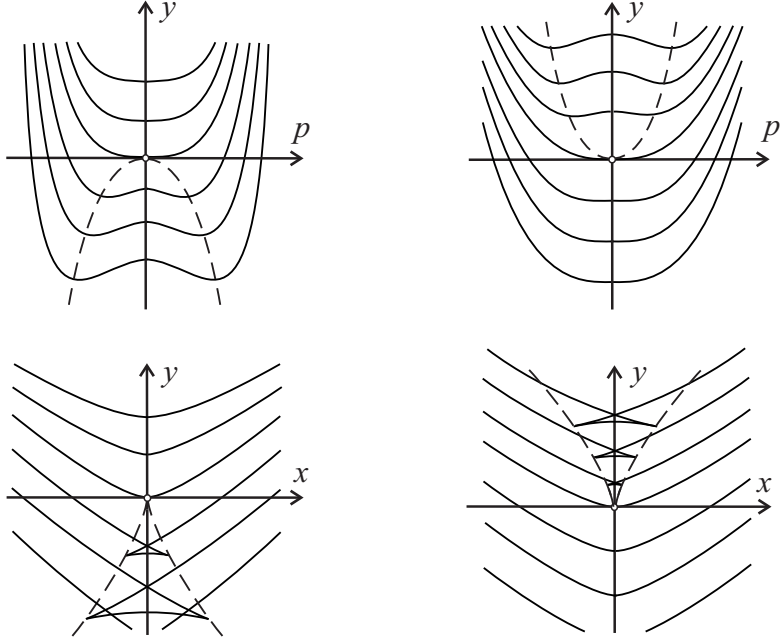}
\caption{
Elliptic pleated point (left) and hyperbolic pleated point (right).
The dashed lines are the criminant and the discriminant curve on the $(x,p)$-plane and the $(x,y)$-plane, respectively.
}\label{Fig0}
\end{center}
\end{figure}

2. \textit{A pleated proper point}: the conditions $F_{pp} = 0$, $F_{ppp} \neq 0$, $F_{xp} \neq 0$ and $F_x \neq 0$ hold true at $\O$. Then the lifted field $X$ is defined, $\KK$ is regular and it has the vertical tangential direction at $\O$, the projection $\pi$ has a pleat at $\O$; see Fig.~\ref{Fig00} (right).
The solutions of equation \eqref{1} in a neighborhood of $\O$ can be described by appropriate sections of the swallow tail \cite{Bruce-84}. There are two essentially different phase portraits in this case, which are called \textit{elliptic pleat} and \textit{hyperbolic pleat}, see Fig.~\ref{Fig0}.
However, there is no a visible classification for equation \eqref{1} in this case, since functional invariants
occur even in the topological category \cite{AI, Dav-85}.

3. \textit{A folded improper point}: the conditions $F_{pp} \neq 0$, $F_y \neq 0$ and $F_x = 0$ hold true at $\O$.
The criminant $\KK$ is regular and not vertical at $\O$, the projection $\pi$ has a fold at all points of $\KK$,
but the lifted field $X$ is not defined (the surface $\FF$ is tangent to the contact plane) at the point~$\O$.

\medskip

Consider the third case in more detail.
Due to $F_y(\O) \neq 0$ equation \eqref{1} can be locally presented in the form
\begin{equation}
\frac{1}{2} \Bigl( ap^2 + 2bxp + cx^2 \Bigr) + f(x,p) = y,   \ \ \ p = \frac{dy}{dx},
\label{3}
\end{equation}
where $a \neq 0$, $b,c$ are real constants, and the germ $f(x,p)$ at $\O$ is 2-flat.\footnote{
The germ of a smooth function is called $k$-flat at $\O$ if its Taylor series at $\O$ starts with monomials of degree greater than $k$.
}
Then formula \eqref{2} for the lifted field on the surface $\FF$ reads
\begin{equation}
\dot x = ap + bx + f_p(x,p), \ \ \
\dot p = (1-b)p - cx - f_x(x,p).
\label{4}
\end{equation}
By $\La$ denote the matrix of the linear part of \eqref{4} at the point $\O$. By $\la_{1,2}$ denote the eigenvalues of $\La$ and by $\ee_{\la_{1,2}}$ denote the corresponding eigendirections (if $\la_{1,2}$ are real and do not coincide). A straightforward computation shows that $\O$ is non-degenerate singular point (saddle, node or focus) of the vector field \eqref{4} if and only if the curves $\KK$ and $\LL$ are regular and transversal at $\O$.

For a generic germ \eqref{3} the following set of conditions holds: both $\la_{1,2}$ are non-zero, $|\la_{1}|:|\la_{2}| \neq 1$, and $\ee_{\la_{1,2}}$ are not tangent neither to $\KK$ nor to the vertical direction at $\O$. (The condition that $\ee_{\la_{1,2}}$ are not tangent to the vertical direction, automatically follows from $a \neq 0$.) Folded improper points satisfying these conditions are called \textit{well-folded}.

In the paper \cite{Dav-85}, A.~Davydov obtained a list of $C^{\infty}$-smooth normal forms of equation \eqref{3}
in a neighborhood of well-folded singular points that satisfy the linearizability condition consisting in that
the germ of the vector field \eqref{4} is $C^{\infty}$-smoothly equivalent to its linear part.
In particular, the linearizability condition holds if between the eigenvalues $\la_{1,2}$ there are no resonance relations
\begin{equation}
\la_i = k_1\la_{1} + k_2\la_{2}, \ \ \ k_{1,2} \in \bZ_+, \ \ k_{1}+k_{2} \ge 2.
\label{5}
\end{equation}
In a neighborhood of a non-resonant well-folded singular point satisfying the linearizability condition,
equation \eqref{3} can be brought to the normal form $(p+\alf x)^2 = y$, where $\alf<0$, $0 < \alf < 1/8$, $\alf > 1/8$ if the point $\O$ is respectively the saddle, the node, the focus of the vector field \eqref{4},
by a $C^{\infty}$-smooth diffeomorphism of the $(x,y)$-plane preserving~$\O$.
Moreover, if we use homeomorphisms of the $(x,y)$-plane, then the parameter $\alf$ can be made to be an arbitrary constant from the corresponding interval.

Normal forms for the resonant singular points (saddles or nodes) were obtained in \cite{Dav-96}.

\medskip

However, if we deal with families of implicit differential equations depending on parameters, some others types of singular points occur generically.

For instance, the case when one of the eigendirections $\ee_{\la_{1,2}}$ is tangent to the criminant (consequently, one of the eigenvalues $\la_{1,2}$ is equal to zero) is considered in \cite{Dav-95}.

Singularities of binary differential equations (which describe the net of principal curvature lines on a surface in Euclidean space) near umbilic points are investigated in \cite{BF-bin, BT-bin, BT-97}.

The case when the surface $\FF$ is not regular at $\O$ (singularity of Morse type) is considered in \cite{BFT-2000}, see also the paper \cite{A-Surfaces}.

A brief survey of these results can be found in \cite{Rem}, see also \cite{Dav-Japan}.

\medskip

In this paper, we investigate the omitted case when the projection $\pi$ has a pleat at a singular point $\O$ of the lifted field. In accordance with our terminology, we call such singular points \textit{pleated improper}.
The surface $\FF$ and the criminant $\KK$ are supposed to be regular, hence the equation is locally equivalent to \eqref{3} with $a=0$, $b \neq 0$.
In this case the eigenvalues $\la_{1,2} = 1-b,b$, and the eigendirection $\ee_{1-b}$ is vertical and tangent to the criminant at~$\O$.

From the aforesaid, it follows that even topological normal forms of implicit differential equations at pleated improper singular points contain
functional invariants, and there is no a visible smooth or topological classification.
However, we prove that there are only six essentially different phase portraits, which are presented in the next section (Fig.~\ref{Fig1} and \ref{Fig2}).

\section{Main results}

Consider implicit differential equation \eqref{3} satisfying the conditions
$$
a=0; \ \ \  b \neq -2, 0, \frac{1}{2}, \frac{2}{3}, 1; \ \ \ f_{ppp}(0,0,0) \neq 0,
$$
in a small neighborhood of $\O$.
The condition $b \neq 0$, $f_{ppp}(0,0,0) \neq 0$ means that the projection $\pi$ has a pleat at $\O$.
The condition $b \neq 0, \frac{1}{2}, 1$ means that the point $\O$ is saddle ($b<0$ or $b>1$) or non-degenerate node ($0<b<1$) of the lifted field $X$.
Finally, $b \neq -2$, $b \neq \frac{2}{3}$ concerns the projection of integral curves of $X$ on the $(x,y)$-plane,
this condition will become clear later on (see Lemma~3).

The five values $b = -2, 0, \frac{1}{2}, \frac{2}{3}, 1$, which we are excluding from consideration, split the range of the parameter $b$ into six intervals corresponding to six different phase portraits of equation \eqref{3}
in a neighborhood of the pleated improper singular point $\O$.

In what follows we use well-known facts from qualitative theory of differential equations, which can be found in \cite{AI}.

\medskip

Without loss of generality we can assume $f_{ppp}(0,0,0) = -2$ (this can be obtained by a scaling of $x$, which does not change neither $a=0$ nor $b$) and $c=0$ (this can be obtained by the change of variables $y \to y - \frac{c}{2(2b-1)}x^2$, which kills the monomial $x^2$). Then the equation reads
\begin{equation}
bxp - \frac{1}{3} p^3 + \phi(p) + x \psi(x,p) = y,   \ \ \ p = \frac{dy}{dx},
\label{6}
\end{equation}
where $\phi(p)$ and $\psi(x,p)$ are 3-flat and 1-flat $C^{\infty}$-germs at $\O$, respectively.

The criminant of equation \eqref{6} is defined by the equality $bx - p^2 + \phi'(p) + x \psi_p(x,p)=0$, which is locally equivalent to $x=p^2 u(p)$ with a $C^{\infty}$-germ $u(p)$ such that $u(0)=b^{-1}$. Substituting this expression in \eqref{6}, we get the asymptotical representation for the criminant and the discriminant curve:
\begin{equation}
x = \frac{1}{b} p^2 + o(p^2),   \ \ \
y = \frac{2}{3} p^3 + o(p^3),  \ \ \ \textrm{as} \ \ p \to 0.
\label{7}
\end{equation}

\medskip

The lifted field $X$ is defined by the vector field
\begin{equation}
\dot x = bx - p^2 + \phi'(p) + x \psi_p(x,p), \ \ \
\dot p = (1-b)p - \psi(x,p) - x\psi_x(x,p),
\label{8}
\end{equation}
which has a saddle or node at $\O$ with the eigenvalues $\la_{1,2}=b,1-b$.
The matrix $\La$ of the linear part of the vector field \eqref{8} at $\O$ is diagonal, and
the eigendirections $\ee_{b}$ and $\ee_{1-b}$ coincide with $\pa_x$ and $\pa_p$, respectively.

In the case of node ($0<b<1$) the resonance relations \eqref{5} have the form $b = \frac{1}{n+1}$, $n \ge 2$ or $b = \frac{n}{n+1}$, $n \ge 3$, for natural $n$. In a neighborhood of the non-resonant and resonant node $\O$, the vector field \eqref{8} is $C^{\infty}$-smoothly orbitally equivalent to
\begin{equation}
\label{100}
\dot \xi  = \xi,  \ \ \ \dot \eta = \bet\eta,  \quad \textrm{where} \quad  \bet =\max \Bigl\{ \frac{b}{1-b}, \frac{1-b}{b} \Bigr\},
\end{equation}
and
\begin{equation}
\label{10}
\dot \xi  = \xi,  \ \ \ \dot \eta = n\eta + \eps \xi^n,  \quad \textrm{where} \quad \eps \in \{0,1\},
\end{equation}
respectively.

\medskip

{\bf Lemma 1.}
{\it
The vector field \eqref{8} has at least one integral curve $\CC \in C^{k}$, $k \ge 2$, passing through $\O$ with the vertical tangential direction $\ee_{1-b}$, and at least one integral curve $\CCC \in C^{l}$, $l \ge 1$, passing through $\O$ with the tangential direction $\ee_{b}$. Moreover,

1) \ $k=l=\infty$ if $\O$ is a saddle or non-resonant node or resonant node with $\eps = 0$,

2) \ $k=\infty$ and $l=n-1$ if $b = \frac{1}{n+1}$, $n \ge 2$, $\eps \neq 0$,

3) \ $k=n-1$ and $l=\infty$ if $b = \frac{n}{n+1}$, $n \ge 3$, $\eps \neq 0$.

}

{\sc Proof.}
In the case of saddle the integral curves $\CC$ and $\CCC$ are separatrices, and the statement is trivial.
Indeed, by the Hadamard--Perron theorem, a $C^{\infty}$-smooth vector field with hyperbolic singular point $\O$ on the plane has $C^{\infty}$-smooth stable and unstable manifolds passing through $\O$ and tangent to $\ee_{\la_{1,2}}$ at this point.

In the case of non-resonant node we have the normal form \eqref{100}, and after integrating obtain the family of integral curves $\eta = c\xi^{\bet}$, $c={\rm const}$, with common tangential direction $\pa_{\xi}$, and the sole integral curve $\xi=0$.
The family $\eta = c\xi^{\bet}$ contains at least one $C^{\infty}$-smooth integral curve (with $c=0$).
Hence the initial vector field has at least one $C^{\infty}$-smooth integral curve tangent to $\ee_{1-b}$ and
at least one $C^{\infty}$-smooth integral curve tangent to $\ee_{b}$.

For the resonant node with $b = \frac{1}{n+1}$, $n \ge 2$, vector field (8) has the normal form \eqref{10},
where the direction $\pa_{\xi}$ corresponds to $\ee_{b}$ and $\pa_{\eta}$ corresponds to $\ee_{1-b}$.
Integrating the differential equation $d\eta/d\xi = (n\eta + \eps \xi^n)/\xi$,
we get the family of integral curves
\begin{equation}
\eta = \xi^n (c + \eps \ln |\xi|), \ \ \  c={\rm const},
\label{11}
\end{equation}
with common tangential direction $\pa_{\xi}$,
and the sole integral curve $\xi=0$.
The integral curves \eqref{11} are $C^{\infty}$-smooth if $\eps = 0$ and $C^{n-1}$-smooth (but not $C^n$-smooth at $\O$) if $\eps \neq 0$.
The integral curve $\xi=0$ corresponds to the integral curve $\CC \in C^{\infty}$, and
any integral curve of the family \eqref{11} corresponds to $\CCC \in C^{l}$ with $l=n-1$ if $\eps \neq 0$ and $l=\infty$ if $\eps = 0$.

For the resonant node with $b = \frac{n}{n+1}$, $n \ge 3$, vector field (8) has the normal form \eqref{10}, where the direction $\pa_{\xi}$ corresponds to $\ee_{1-b}$ and $\pa_{\eta}$ corresponds to $\ee_{b}$. The integral curve $\xi=0$ corresponds to the integral curve $\CCC \in C^{\infty}$, and any integral curve of the family \eqref{11} corresponds to $\CC \in C^{k}$ with $k=n-1$ if $\eps \neq 0$ and $k=\infty$ if $\eps = 0$.
\ebox

\medskip

{\bf Lemma 2.}
{\it
Let $\CC$ be a $C^2$-smooth integral curve of the vector field \eqref{8} passing through the singular point $\O$ with the vertical tangential direction $\ee_{1-b}$.
Then the germ of $\CC$ is given by $x = v_0 p^2 + o(p^2)$, where $v_0={1}/{(3b-2)}$,
and the projection $\pi(\CC)$ has the asymptotical representation
\begin{equation}
x = \frac{1}{3b-2} p^2 + o(p^2),   \ \ \
y = \frac{2}{3(3b-2)} p^3 + o(p^3),  \ \ \ \textrm{as} \ \ p \to 0.
\label{9}
\end{equation}
}

{\sc Proof.}
The integral curve $\CC$ can be locally presented in the form $x = v_0 p^2 + v(p)$ with a $C^2$-smooth germ $v(p)=o(p^2)$ as $p \to 0$.
Substituting this expression in the equality
$$
\frac{dx}{dp} = \frac{bx - p^2 + \phi'(p) + x \psi_p(x,p)}{(1-b)p - \psi(x,p) - x\psi_x(x,p)},
$$
we get
$$
2v_0p + o(p) = \frac{bv_0p^2 - p^2 + o(p^2)}{(1-b)p + o(p)} \ \ \Rightarrow \ \
 2v_0 = \frac{bv_0 - 1}{1-b} \ \ \Rightarrow \ \ v_0=\frac{1}{3b-2}.
$$
Substituting $x = v_0 p^2 + v(p)$ with $v_0={1}/{(3b-2)}$ in \eqref{6}, we get the representation \eqref{9}.
\ebox

\medskip

{\bf Theorem.}
{\it
In a neighborhood of the point $\O$, equation \eqref{6} can be reduced to the form
\begin{equation}
p\, \Bigl( bx - \frac{1}{3} p^2 + \tphi(p) + x \tpsi(x,p) \Bigr) = y,   \ \ \ p = \frac{dy}{dx},
\label{12}
\end{equation}
where the germ $\tphi(p)$ is 2-flat at $\O$ and $\tpsi(0,0)=0$, by appropriate change of variables
\begin{equation}
y \to y - u(x), \ \ \ u(0)=u'(0)=0, \ \ \ u(\cdot) \in C^s,
\label{13}
\end{equation}
where $s=n$ if $b = \frac{1}{n+1}$, $n \ge 2$, $\eps \neq 0$, and $s={\infty}$ otherwise.
}

{\sc Proof.}
The statement is equivalent to the existence of a solution $y = u(x)$ such that $u(\cdot) \in C^s$ and $u(0)=u'(0)=0$. Indeed, the change of variables \eqref{13} takes the solution $y = u(x)$ to $y=0$, consequently, it brings equation \eqref{6} to the form \eqref{12}.

Thus it is necessary and sufficient to establish the existence of a $C^{s-1}$-smooth integral curve of the vector field \eqref{8} passing through the point $\O$ with the tangential direction $\ee_b$. Clearly, the integral curve $\CCC$ from Lemma~1 satisfies the required conditions.
\ebox

\medskip

The representations \eqref{7} and \eqref{9} show that the curves $\pi(\KK)$ and $\pi(\CC)$ are semicubic parabolas on the $(x,y)$-plane having the common cusp at $\O$ with the same tangential direction $\pa_x$. To determine the mutual arrangement of $\pi(\KK)$ and $\pi(\CC)$ in a neighborhood of $\O$, it is convenient to represent the semicubic parabolas $\pi(\KK)$ and $\pi(\CC)$ in the form of a sole algebraic equation:
\begin{equation}
x = A p^2, \ \ y = B p^3 \ \ \Rightarrow \ \ y^2 = m x^3, \ \ m =\frac{B^2}{A^3}.
\label{2013}
\end{equation}

\medskip

{\bf Lemma 3.}
{\it
The semicubic parabolas $\pi(\CC)$ and $\pi(\KK)$ belong to different semiplanes into which the $y$-axis divides the $(x,y)$-plane if $0 < b < \frac{2}{3}$, and to the same semiplane otherwise.
Moreover, $\pi(\CC)$ lies in the ``smaller'' (tongue-like) of the domains into which $\pi(\KK)$ locally divides the $(x,y)$-plane if $b<-2$ or $\frac{2}{3}< b$, and vice versa if $-2<b<0$.
}

{\sc Proof.}
Formula \eqref{2013} gives $m = \frac{4}{9}b^3$ and $m = \frac{4}{9}(3b-2)$ for the semicubic parabolas
defined by the principal parts of asymptotic formulae \eqref{7} and \eqref{9}, respectively.
\ebox

\medskip

\begin{figure}[h!]
\begin{center}
\includegraphics{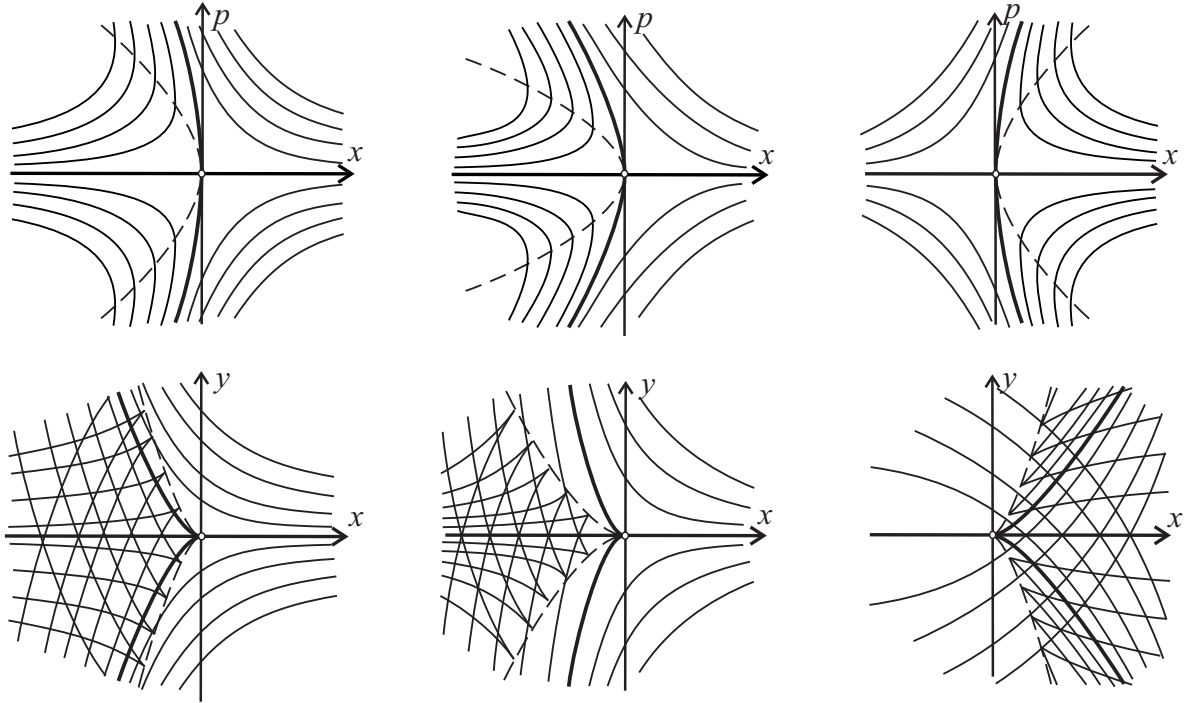}
\caption{
Phase portraits in the cases $S_1, S_2, S_3$ form left to right.
Bold lines are the integral curve $\CC$ and its projection $\pi(\CC)$,
the dashed lines are the criminant and the discriminant curve on the $(x,p)$-plane and the $(x,y)$-plane, respectively.
}\label{Fig1}
\end{center}
\end{figure}

\begin{figure}[h!]
\begin{center}
\includegraphics{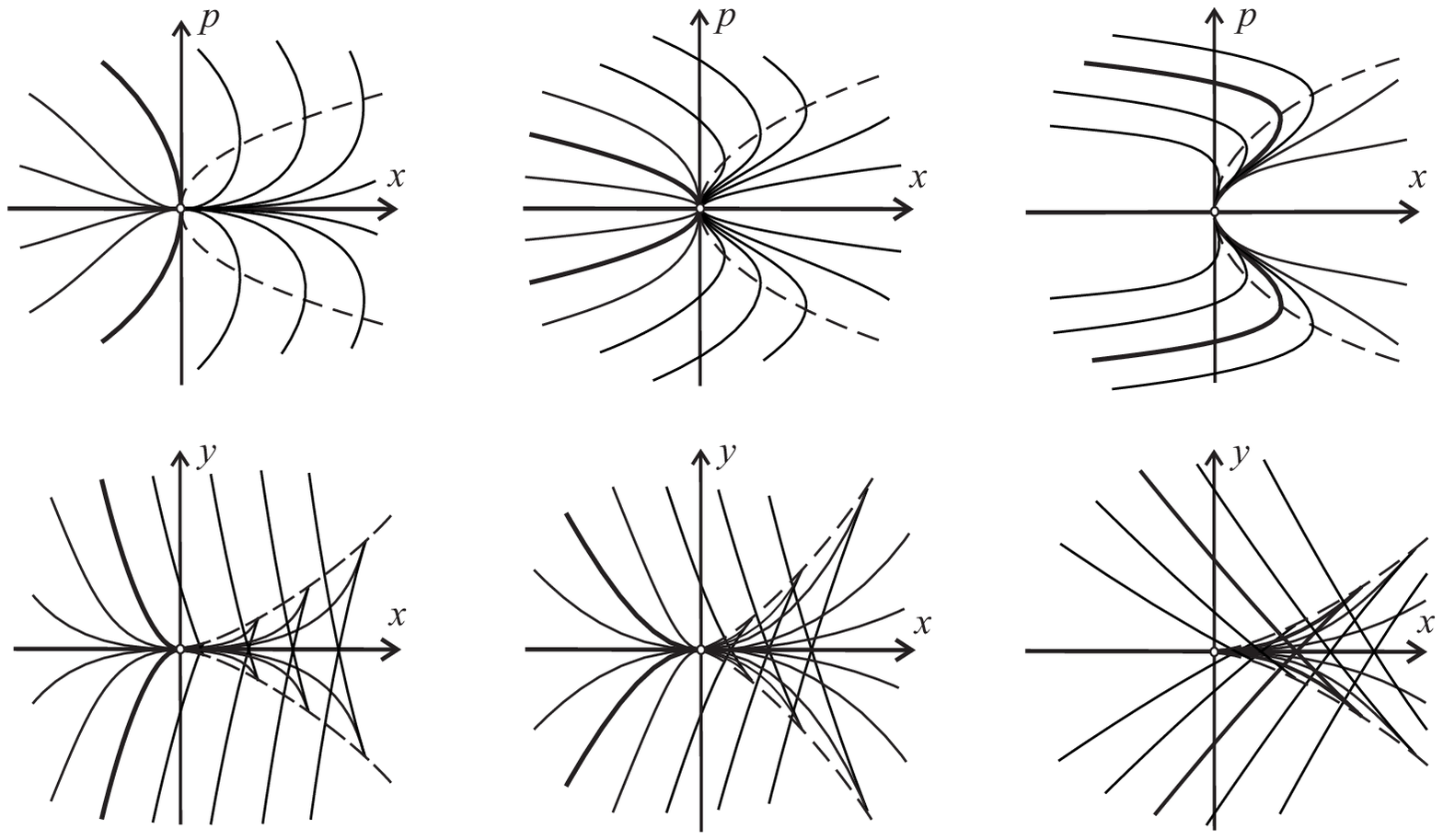}
\caption{
Phase portraits in the cases $N_1, N_2, N_3$ form left to right.
Bold lines are the integral curve $\CC$ and its projection $\pi(\CC)$,
the dashed lines are the criminant and the discriminant curve on the $(x,p)$-plane and the $(x,y)$-plane, respectively.
}\label{Fig2}
\end{center}
\end{figure}

In accord with the above reasonings, classification of phase portraits of equation \eqref{6} in a neighborhood of the pleated improper point $\O$ is presented in Tab.~1 and Fig.~\ref{Fig1}, \ref{Fig2}.

\begin{table}[htb]
\begin{center}
\begin{tabular}{|c|c|c|c|c|c|c|}
\hline
{Name of case} & $S_1$ & $S_2$ & $N_1$ & $N_2$ & $N_3$ & $S_3$ \\
\hline
{Range} & $b<-2$ & $-2<b<0$ & $0<b<\frac{1}{2}$ & $\frac{1}{2}<b<\frac{2}{3}$ & $\frac{2}{3}<b<1$ & $1<b$ \\
\hline
Lifted field & saddle & saddle & node & node & node & saddle \\
\hline
$\sign(\frac{1}{b}):\sign(\frac{1}{3b-2})$ & $-:-$ & $-:-$ & $+:-$ & $+:-$ & $+:+$ & $+:+$ \\
\hline
$|\frac{1}{b}|:|\frac{1}{3b-2}|$ & $>1$ & $>1$ &    &    & $<1$ & $>1$ \\
\hline
$|b^3|:|3b-2|$ & $>1$ & $<1$ &    &    & $>1$ & $>1$ \\
\hline
\end{tabular}
\end{center}
\caption{Classification of phase portraits in a neighborhood of a pleated improper singular point.}
\end{table}

The difference between the cases $N_2$ and $N_3$ needs to be commented. In both cases $b>1-b$, and almost all integral curves of the vector field \eqref{8} have vertical tangential direction $\ee_{1-b}$ at $\O$. Resonance relations do not occur in the case $N_2$. In the case $N_3$ they have the form $b = \frac{n}{n+1}$, $n \ge 3$, and all integral curves \eqref{11} are at least $C^{2}$-smooth.
In the non-resonant case, the germ of the vector field \eqref{8} has the orbital normal form \eqref{100}
with $1<\bet<2$ if $\frac{1}{2}<b<\frac{2}{3}$ and $\bet>2$ if $\frac{2}{3}<b<1$.
Hence the curves $\eta = c\xi^{\bet}$, $c \neq 0$, are $C^2$-smooth in the case $N_3$ and only $C^1$-smooth in the case $N_2$.
This also holds true for the corresponding integral curves of the vector field \eqref{8}.

In the case $N_3$ Lemma~2 is applicable to all integral curves of \eqref{8} passing through the point $\O$ with vertical tangential direction $\ee_{1-b}$, and all the curves have the same 2-jet at $\O$.
In the case $N_2$ Lemma~2 is not applicable to all integral curves of \eqref{8} passing through $\O$ with vertical tangential direction $\ee_{1-b}$
(except for only one curve mentioned in Lemma~1), and this family contains both convex and concave curves.

The phase portraits of equation \eqref{12} for all cases from the table are presented below,
where the integral curve $\CCC$ from Lemma~1 coincides with the axis $p=0$, and its projection $\pi(\CCC)$ coincides with the axis $y=0$.

\section*{Acknowledgements}
A.O. Remizov was supported by grant from FAPESP, proc. 2012/03960-2 for visiting ICMC-USP, S\~ao Carlos (Brazil).
He expresses deep gratitude to prof. Farid Tari for hospitality and useful discussions.


\end{document}